\documentclass[11pt,leqno]{amsart}
\usepackage{epsfig}
\usepackage{amssymb}
\usepackage{amscd}
\usepackage[matrix,arrow]{xy}

\newtheorem{theo}{Theorem}[section]
\newtheorem{lemma}[theo]{Lemma}
\newtheorem{defi}[theo]{Definition}
\newtheorem{prop}[theo]{Proposition}

\newtheorem{cor}[theo]{Corollary}
\newtheorem{remark}[theo]{Remark}

\newtheorem{example}[theo]{Example}
\numberwithin{equation}{section}

\def\pre-tr{\operatorname{pre-tr}}

\def\Hom{\operatorname{Hom}}

\newcommand{\bbZ}{{\mathbb Z}}

\newcommand{\bbQ}{{\mathbb Q}}

\newcommand{\cO}{{\mathcal O}}

\newcommand{\cA}{{\mathcal A}}
\newcommand{\cB}{{\mathcal B}}

\newcommand{\cC}{{\mathcal C}}

\newcommand{\cH}{{\mathcal H}}

\newcommand{\Ker}{\operatorname{Ker}}
\newcommand{\im}{\operatorname{Im}}

\newcommand{\id}{\operatorname{id}}

\newcommand{\Coder}{\operatorname{Coder}}
\newcommand{\pd}{\operatorname{pd}}
\newcommand{\depth}{\operatorname{depth}}
\newcommand{\Der}{\operatorname{Der}}
\newcommand{\ad}{\operatorname{ad}}

\title[Formality of DG algebras (after Kaledin)]{Formality of DG algebras (after Kaledin)}

\author{Valery A.~Lunts}
\address{Department of Mathematics, Indiana University,
Bloomington, IN 47405, USA} \email{vlunts@indiana.edu}

\thanks{The author was partially supported by the NSA grant H98230-07-1-0071}

\begin{document}

\begin{abstract} We provide proper foundations and proofs for the
main results of [Ka]. The results
 include a flat base change for
 formality and behavior of formality in flat families of $A(\infty)$ and DG algebras.
\end{abstract}

\maketitle

\section{introduction}

Let $k$ be a field of characteristic zero. Given a DG algebra $\cA$
over $k$ Kaledin [Ka] defines a cohomology class $K_A$ which
vanishes if and only if $\cA$ is formal. (This class $K_A$ is an
element of the second Hochschild cohomology group of a DG algebra
$\tilde{\cA}$ which is closely related to $\cA$.) This is a
beautiful result which has many important applications. One of the
applications is mentioned in [Ka] (Theorem 4.3): if one has a "flat"
family $\cA _X$ of DG algebras parametrized by a scheme $X$, then
formality of the fiber $\cA _x$ is a closed condition on $x\in X$.

Unfortunately, the paper [Ka] is hard to read. There are many
misprints and inaccuracies. The definition and treatment of the
Hochschild cohomology of a family of DG algebras is unsatisfactory:
for example, in the proof of main Theorem 4.3 it is implicitly
assumed that the Hochschild cohomology behaves well with respect to
specialization.

But nonetheless we found the paper [Ka] inspiring and decided to
provide the necessary foundations and proofs of its main results.

Unlike [Ka] we found it more convenient to work with $A(\infty)$
algebras rather than with DG algebras. Namely, for a commutative
ring $R$ we consider  $A(\infty)$ $R$-algebras which are minimal
($m_1=0$) and {\it flat}, i.e. each $R$-module $H^n(A)=A^n$ is
projective. That is what we mean by a {\it flat family} of
$A(\infty)$ algebras over $SpecR$. We are mostly interested in the
case when the $R$-module $A$ is finite.

The behavior of the ($R$-linear) Hochschild cohomology $HH_R(A)$
with respect to base change $R\to Q$ is hard to control. For
$A(\infty)$ algebras $A$ which are {\it finitely defined} (i.e. only
finitely many $m_i$'s are not zero) one may consider the Hochschild
cohomology with {\it compact supports} $HH_{R,c}(A)$. It comes with
a natural map $HH_{R,c}(A)\to HH_R(A)$ which is injective in cases
which are important for us. The groups $HH_{R,c}$ have better
behavior with respect to base change and they contain Kaledin's
cohomology classes, which are obstructions to formality. Thus in
essential places we work with $HH_c(A)$ and not with $HH(A)$. The
good functorial behavior of $HH_c(A)$ allows us to prove a
faithfully flat base change result for formality (Proposition 6.2).
A similar result for {\it commutative} DG algebras over a field was
proved by Sullivan [Su] (see also [HaSt]).

The paper is organized as follows. In Section 2 we recall
$A(\infty)$ algebras over arbitrary commutative rings, their bar
constructions, quasi-isomorphisms and Kadeishvili's theorem. We also
relate the DG formality of flat DG algebras to $A(\infty)$ formality
of their minimal models. In Section 3 we recall Hochschild
cohomology, introduce Hochschild cohomology with compact supports
and discuss its properties.  In Section 4 we define Kaledin's
cohomology class and discuss its relation to (infinitesimal)
formality. In Section 5 we consider the "deformation to the normal
cone" $\tilde{A}$ of an $A(\infty)$ algebra $A$ and prove the
Kaledin's key result. Section 6 contains applications of this result
to the behavior of formality in flat families of $A(\infty)$ (or DG)
algebras. Finally in Section 7 we define Kaledin cohomology class in
the general context of DG Lie algebras.

We thank Dima Kaledin for answering many questions about [Ka] and
Bernhard Keller and Michael Mandell for answering general questions
about $A(\infty)$ algebras. Jee Koh helped us with commutative
algebra. We should mention the paper [Hi] by Vladimir Hinich which
helped us understand what Kaledin was trying to do. We also thank
the anonymous referee for several useful remarks and suggestions.

\section{$A(\infty)$ algebras}

A good introduction to $A(\infty)$ algebras is [Ke]. However there
seems to be no systematic treatment of $A(\infty)$ algebras over an
arbitrary commutative ring.

\subsection{$A(\infty)$-algebras} Fix a commutative unital ring $R$. The sign $\otimes $ means $\otimes _R$.
 We want to study $A(\infty)$ $R$-algebras and quasi-isomorphisms between them. Let us
recall the definitions.

Let $A=\oplus _{n\in \bbZ}A^n$ be a graded $R$-module. A structure
of an $A(\infty)$ $R$-algebra (or, simply, $A(\infty)$ algebra) on
$A$ is a collection $m=(m_1,m_2,...)$, where $m_i:A^{\otimes i}\to
A$ is a homogeneous $R$-linear map of degree $2-i$. The maps
$\{m_i\}$ must satisfy for each $n\geq 1$ the following identity:
 $$\sum (-1)^{r+st}m_u(1^{\otimes r}\otimes m_s\otimes 1^{\otimes t})=0,$$
where the sum runs over all decompositions $n=r+s+t$ and we put
$u=r+1+t$.

We denote the resulting $A(\infty)$-algebra by $(A,m),$ or
$(A,(m_1,m_2,...))$ or simply by $A$.

\begin{itemize}
\item If $m_i=0$ for $i\neq 2$, then $A$ is simply a graded
associative $R$-algebra.
\item If $m_i=0$ for $i\neq 1,2$ then $A$ is a DG $R$-algebra.
\item If $m_1=0$ then $A$ is called {\it minimal}. Note that in this
case $A$ is in particular a graded associative $R$-algebra with
multiplication $m_2$.
\item In any case $A$ is a complex of $R$-modules with the
differential $m_1$ and the cohomology $H(A)$ is a graded associative
$R$-algebra with multiplication defined by $m_2$.
\end{itemize}

\subsection{$A(\infty)$-morphisms} Given $A(\infty)$ algebras $A,B$ an $A(\infty)$ {\it
morphism} $f:A\to B$ is a collection $f=(f_1,f_2,...)$, where
$f_i:A^{\otimes i}\to B$ is an $R$-linear map of degree $1-i$ such
that for each $n\geq 1$ the following identity holds.
$$\sum (-1)^{r+st}f_u(1^{\otimes r}\otimes m_s\otimes 1^{\otimes
t})=\sum (-1)^sm_r(f_{i_1}\otimes f_{i_2}\otimes ...\otimes
f_{i_r}),$$ where the first sum runs over all decompositions
$n=r+s+t$, we put $u=r+1+t$, and the second sum runs over all $1\leq
r\leq n$ and all decompositions $n=i_1+...+i_r$; the sign on the
right hand side is given by
$$s=(r-1)(i_1-1)+(r-2)(i_2-1)+...+2(i_{r-2}-1)+(i_{r-1}-1).$$

\begin{itemize}
\item We have $f_1m_1=m_1f_1$, i.e. $f_1$ is a morphism of
complexes.
\item We have
$$f_1m_2=m_2(f_1\otimes f_1)+m_1f_2+f_2(m_1\otimes 1+1\otimes
m_1),$$ which means that $f_1$ commutes with the multiplication
$m_2$ up to a homotopy given by $f_2$. In particular, if $A$ and $B$
are minimal, then $f_1$ is a homomorphism of associative algebras
$f_1:(A,m_2)\to (B,m_2)$.
\end{itemize}

We call $f$ a {\it quasi-isomorphism} if $f_1:A\to B$ is a
quasi-isomorphisms of complexes. $f$ is called the identity
morphism, denoted $\id$, if $A=B$ and $f=(f_1=\id,0,0,...)$.

Let $C$ be another $A(\infty)$ algebra and $g=(g_1,g_2,...):B\to C$
be an $A(\infty)$ morphism. The composition $h=g\cdot f:A\to C$ is
an $A(\infty)$-morphism which is defined by
$$h_n =\sum (-1)^sf_r(g_{i_1}\otimes ...\otimes g_{i_r}),$$
where the sum and the sign are as in the defining identity.

$A(\infty)$ algebras $A$ and $B$ are called quasi-isomorphic if
there exists $A(\infty)$ algebras $A_1, A_2, ..., A_n$ and
quasi-isomorphisms
$$A\leftarrow A_1\rightarrow ...\leftarrow
A_{n}\rightarrow B.$$

An $A(\infty)$ algebra $A$ is called {\it formal} if it is
quasi-isomorphic to the $A(\infty)$ algebra $(H(A),(0,m_2,0,...))$.

\subsection{Bar construction}
The notions of $A(\infty)$ algebra and $A(\infty)$ morphism can be
compactly and conveniently described in terms of the bar
construction.

 Let $A$ be a graded $R$-module, $A[1]$
 its shift $A[1]^n=A^{n+1}$. Let
$$\overline{T}A[1]=\bigoplus_{i\geq 1}A[1]^{\otimes i}$$
be the reduced cofree $R$-coalgebra on the $R$-module $A[1]$ with
the comultiplication
$$\Delta (a_1,...,a_n)=\sum_{i=1}^{n-1}(a_1,...,a_i)\otimes
(a_{i+1},...,a_n)$$ so that $\Delta (a)=0$ and $\Delta (a_1, a_2)=a_1\otimes a_2.$
 Denote by $\Coder(\overline{T}A[1])$ the
graded $R$-module of homogeneous $R$-linear coderivations of the
coalgebra $\overline{T}A[1]$. The composition of a coderivation with
the projection to $T^1A[1]=A[1]$ defines an isomorphism of graded
$R$-modules
$$\Coder (\overline{T}A[1])\simeq \Hom _R(\overline{T}A[1],A[1]).$$

Thus a coderivation of degree $p$ is determined by a collection
$(d_1,d_2,...)$, where $d_i:A[1]^{\otimes i}\to A[1]$ is an
$R$-linear map of degree $p$.

Denote by $s:A\to A[1]$ the shift operator. Given an $R$-linear map
$m_i:A^{\otimes i}\to A$ of degree $2-i$ we can define an $R$-linear
map $d_i:A[1]^{\otimes i}\to A[1]$ of degree 1 by commutativity of
the following diagram
$$\begin{array}{ccc}
A^{\otimes i} & \stackrel{m_i}{\to} & A\\
\downarrow s^{\otimes i} & & \downarrow s\\
A[1]^{\otimes i} & \stackrel{d_i}{\to} & A[1]
\end{array}
$$
Thus $d_i(sa_1\otimes ...\otimes sa_i)=(-1)^nsm_i(a_1\otimes
...\otimes a_i),$ where
$n=\frac{i(i+1)}{2}+(i-1)\deg(a_1)+(i-2)\deg(a_2)+...+\deg(a_{i-1}).$
Then a collection of $R$-linear maps $m=(m_1,m_2,...),$
$m_i:A^{\otimes i}\to A$ of degree $2-i,$ defines a structure of an
$A(\infty)$ $R$-algebra on $A$ if and only if the corresponding
collection $d=(d_1,d_2,...),$ $d_i:A[1]^{\otimes i} \to A[1]$ of
degree $1,$ defines an $R$-linear coderivation of the coalgebra
$\overline{T}A[1]$ such that $d^2=0.$ Given an $A(\infty)$ algebra
$(A,m)$ we will also denote (abusing notation) by the same letter
$m$ the corresponding coderivation of $\overline{T}A[1].$ The
resulting DG coalgebra $(\overline{T}A[1],m)$ is called the {\it bar
construction} of $A$ and is denoted $\cB A.$

Let $B$ be another $A(\infty)$ algebra. In a similar manner (using
appropriate sign changes) there is a bijection between $A(\infty)$
morphisms $A\to B$ and homomorphisms of degree zero of DG coalgebras
$\cB A\to \cB B.$ Again we will usually use the same notation for
both.

Let $f=(f_1,f_2,...):\overline{T}A[1]\to \overline{T}B[1]$ be a
homomorphism of coalgebras. Then for each $n$
$$f(\bigoplus _{i\leq n}T^iA[1])\subset \bigoplus _{i\leq n}T^iB[1].$$
The map $f$  is an isomorphism if and only if $f_1$ is an
isomorphism. On the other hand if $f_1=0$ and $A=B$, then the map
$f$ is locally nilpotent.

Similar considerations apply to coderivations
$g=(g_1,g_2,...):\overline{T}A[1]\to \overline{T}A[1]$. Namely, let
$g$ have degree zero and $g_1=0$, then $g$ is locally nilpotent and
hence the coalgebra automorphism
$$\exp(g):\overline{T}A[1]\to \overline{T}A[1]$$
is well defined (provided $\bbQ \subset R$).

\subsection{Flat $A(\infty)$ algebras and their minimal models}

\begin{defi}
An $A(\infty)$ $R$-algebra $A$ is called {\it flat}  if each
cohomology $H^i(A)$ is a projective $R$-module.
\end{defi}

Thus if $R$ is a field then any $A(\infty) $ algebra if flat. We
consider a flat $A(\infty)$ $R$-algebra as a flat family of
$A(\infty)$ algebras over $SpecR$. Let us recall the following
simple important result of Kadeishvili.

\begin{theo} [Kad1] Let $A$ be a flat $A(\infty)$ $R$-algebra.
Choose a quasi-isomorphism of complexes of $R$-modules $g:H(A)\to A$
(the differential in $H(A)$ is zero). Then there exists a structure
of a minimal $A(\infty)$ algebra on $H(A)$ with $m_2$ being induced
by the $m_2$ of $A$ and an $A(\infty)$ morphism $f=(g=f_1,f_2,...)$
from $H(A)$ to $A$ (which is a quasi-isomorphism).
\end{theo}

We call the $A(\infty)$ algebra $H(A)$ as in the above theorem a
minimal model of $A.$

Let $A$ and $B$ be $A(\infty)$ $R$-algebras and $f,g$ morphisms from
$A$ to $B.$ Let $F, G$ denote the corresponding morphisms of DG
coalgebras $\cB A\to \cB B.$ One defines $f$ and $g$ to be homotopic
if $F$ and $G$ are homotopic, i.e. if there exists a homogeneous
$R$-linear map $H:\cB A\to \cB B$ of degree $-1$ such that
$$\Delta \cdot H=F\otimes H +H\otimes G\quad \text{and} \quad
F-G=m_B\cdot H +H \cdot m_A.$$

\begin{lemma} In the above notation assume that $A$ and $B$ are
minimal. Let $f:A\to B$ and $g:B\to A$ be morphisms such that
$g\cdot f$ and $f\cdot g$ are homotopic to the identity (i.e. $A$
and $B$ are homotopy equivalent). Then the corresponding morphisms
$F:\cB A \to \cB B$ $G:\cB B\to \cB A$ are mutually inverse
isomorphisms.
\end{lemma}

\begin{proof} Let $H:\cB A \to \cB A [-1]$ be a homotopy between
morphisms $G\cdot F$ and $\id _{\cB A}.$ Then $H$ is defined by a
collection of $R$-linear maps $h_i:A[1]^{\otimes i}\to A[1],$ $i\geq
1,$ which satisfy some properties ([Le-Ha],1.2.1.7).

Let $F=(f_1,f_2,...),\ G=(g_1,g_2,...),\ G\cdot F =(t_1,t_2,...).$
Then
$$(G\cdot F)\vert _{A[1]}=t_1\vert _{A[1]}=g_1\cdot f_1 \vert
_{A[1]}.$$ Also $H\vert _{A[1]}=h_1\vert _{A[1]}.$ Since $A$ is
minimal the equation
$$G\cdot F -\id =m \cdot H+H\cdot m$$
when restricted to $A[1]$ becomes $g_1\cdot f_1 -\id =0\cdot h_1
+h_1\cdot 0.$ So $t_1=\id,$ i.e. $G\cdot F:\cB A\to \cB A$ is an
automorphism.
\end{proof}

Let us recall another result of Kadeishvili.

\begin{theo} [Kad2] a) Homotopy is an equivalence relation on the
set of morphisms of $A(\infty)$ algebras $A\to B.$

Denote by $\cH$ the category obtained by dividing the category of
$A(\infty)$ $R$-algebras by the homotopy relation.

b) Assume that $C$ is an $A(\infty)$ $R$-algebra such that the
$R$-module $C^n$ is projective for all $n\in \bbZ.$ Then a
quasi-isomorphism of $A(\infty)$ algebras $s:A\to B$ induces an
isomorphism
$$s_*:Hom _{\cH}(C,A)\to \Hom _{\cH}(C,B).$$
\end{theo}

\begin{cor} On the full subcategory of $A(\infty)$ $R$-algebras
which consists of algebras $C$ such that the $R$-module $C^n$ is
projective for all $n\in \bbZ$ the relation of quasi-isomorphism
coincides with the relation of homotopy equivalence.
\end{cor}

\begin{cor} Let $A$ and $B$ be two minimal flat $A(\infty)$
$R$-algebras. Then they are quasi-isomorphic if and only if their
bar constructions $\cB A$ and $\cB B$ are isomorphic.
\end{cor}

\begin{proof} The "if" direction is obvious.

Assume that $A$ and $B$ are quasi-isomorphic, i.e. there exists a
chain of morphisms of $A(\infty)$ $R$-algebras which are
quasi-isomorphisms:
$$A\stackrel{f}{\leftarrow} A_1 \stackrel{g}{\to}A_2 \leftarrow
...\to B.$$ Choose a flat minimal $A(\infty)$ $R$-algebra
$A_1^\prime$ and a quasi-isomorphism $i:A^\prime _1\to A_1.$ The
quasi-isomorphism $f\cdot i:A^\prime _1\to A$ between two minimal
flat $A(\infty)$ algebras induces an isomorphism of their bar
constructions $\cB A^\prime _1\to \cB A.$

Choose a flat minimal $A(\infty)$ algebra $A^\prime _2$ and a
quasi-isomorphism $j:A^\prime _2 \to A_2.$ It follows from Theorem
2.4 b) that the induced maps
$$(g\cdot i)_*:\Hom _{\cH}(C,A^\prime _1)\to \Hom
_{\cH}(C,A_2)\leftarrow \Hom _{\cH}(C,A^\prime _2):j_*$$ are
isomorphisms if $C$ is a minimal flat $A(\infty)$ algebra. In
particular $A^\prime _1$ and $A^\prime _2$ are homotopy equivalent,
hence $\cB A^\prime _1 \simeq \cB A^\prime _2$ by Lemma 2.3
Continuing this way we arrive at an isomorphism $\cB A \simeq \cB
B.$
\end{proof}

This last corollary implies in particular that for a flat
$A(\infty)$ algebra its minimal model (as in Theorem 2.2) is unique
up to a quasi-isomorphism. We will identify a quasi-isomorphism
$A\stackrel{\sim}{\to }B$ between flat minimal $A(\infty)$ algebras
with the corresponding isomorphism $\cB A\stackrel{\sim}{\to}\cB B.$

\begin{cor} Let $A$ be a flat $A(\infty)$ algebra and $B$ be
a minimal flat $A(\infty)$ algebra which is quasi-isomorphic to $A.$
Then there exists a morphism $B\to A$ which is a quasi-isomorphism.
\end{cor}

\begin{proof} Let $H(A)$ be a minimal flat $A(\infty)$ algebra with
a quasi-isomorphism $H(A)\to A$ as in Theorem 2.2. Then the minimal
flat $A(\infty)$ algebras $B$ and $H(A)$ are quasi-isomorphic. So by
Corollary 2.6 $\cB B\simeq \cB H(A).$ Hence there also exists a
quasi-isomorphism $B \stackrel{\sim}{\to}A.$
\end{proof}

\subsection{Flat DG algebras} A DG ($R$-)algebra is an $A(\infty)$
algebra $(A,m)$ such that $m_i=0$ for $i>2.$ A morphism of DG
algebras is a homomorphism of graded associative algebras which
commutes with the differentials. Thus the category of DG algebras is
not a full subcategory of $A(\infty)$ algebras. We say that DG
algebras $A$ and $B$ are DG quasi-isomorphic if there exists a chain
of morphisms of DG algebras
$$A\leftarrow A_1 \to ...\leftarrow A_n \to B$$
where all arrows are quasi-isomorphisms. It is well known that if
$R$ is a field then two DG algebras are quasi-isomorphic (as
$A(\infty)$ algebras) if and only if they are DG quasi-isomorphic.
For a general ring $R$ we have a similar result for flat DG algebras
 (Definition 2.1).

\begin{prop} Let $E$ and $F$ be flat DG $R$-algebras and $A$ and $B$
be their minimal $A(\infty)$-models (Theorem 2.2). The following
assertions are equivalent.

a) $E$ and $F$ are DG quasi-isomorphic.

b) $E$ and $F$ are quasi-isomorphic.

c) $A$ and $B$ are quasi-isomorphic.

d) $\cB A$ and $\cB B$ are isomorphic.
\end{prop}

\begin{proof} Clearly a) $\Rightarrow$ b) and by definition b) $\Leftrightarrow$ c).  Corollary 2.6 implies
that c) $\Leftrightarrow$ d). So it remains to prove that d)
$\Rightarrow $ a).

So assume that $\cB A\simeq \cB B.$  Choose a DG algebra $\tilde{E}$
such that the $R$-module $\tilde{E}^n$ is projective for all $n\in
\bbZ,$ and a DG quasi-isomorphism $\tilde{E}\to E$ (for example
$\tilde{E}$ may be a cofibrant replacement of $E$).

By Corollary 2.7 there exists an $A(\infty)$ morphism $A\to
\tilde{E}$ which is a quasi-isomorphism. By Corollary 2.5 this is a
homotopy equivalence, i.e. the induced morphism of the bar
constructions $\cB A \to \cB \tilde{E}$ is a homotopy equivalence.

Consider the cobar construction $\Omega $ which is a functor from DG
coalgebras to DG algebras [Le-Ha],1.2.2. It is the left adjoint to
the bar construction $\cB.$ The same proof as of Lemma 1.3.2.3 in
[Le-Ha] shows that the adjunction morphism of DG algebras $\Omega
\cB \tilde{E}\to \tilde{E}$ is a quasi-isomorphism.

But the DG coalgebras $\cB A $ and $\cB \tilde{E}$ are homotopy
equivalent. Hence their cobar constructions are also homotopy
equivalent and in particular the DG algebras $\Omega \cB A$ and
$\Omega \cB \tilde{E}$ are DG quasi-isomorphic. (The notion of
homotopy between morphisms of DG algebras is defined for example in
[Le-Ha],1.1.2.) Thus the DG algebras $\Omega \cB A$ and $E$ are DG
quasi-isomorphic.

Similarly one shows that the DG algebras $\Omega \cB B$ and $F$ are
DG quasi-isomorphic. But an isomorphism of DG coalgebras  $\cB A
\simeq \cB B$ induces an isomorphism of DG algebras $\Omega \cB
A\simeq \Omega \cB B.$ This proves the proposition.
\end{proof}

A DG algebra in called DG formal if it is DG quasi-isomorphic to a
DG algebra with the zero differential.

\begin{cor} Let $E$ be a flat DG $R$-algebra with a minimal
$A(\infty)$ model $A.$ Then $E$ is DG formal if and only in $A$ is
formal (Subsection 2.2). So $E$ is DG formal if and only if it is
$A(\infty)$ formal.
\end{cor}

\begin{proof} This follows from the equivalence of a) and c) in
Proposition 2.8.
\end{proof}

In what follows we will be interested only in flat $A(\infty)$ or DG
algebras and hence will usually work with their minimal models.

\section{Hochschild cohomology}

We assume that $A$ is a minimal flat $A(\infty)$ $R$-algebra.

\subsection{}
Consider the graded $R$-module $\Coder (\overline{T}A[1])$ with the
self map of degree 1 given by $d\mapsto [m_A,d]=m_A\cdot
d-(-1)^{\deg d}d\cdot m_A$. Since $m_A^2=0$ this makes $\Coder
(\overline{T}A[1])$ a complex of $R$-modules which we denote by
$C^\bullet _R(A)$.  This complex is called the {\it Hochschild
complex} of $A$. Its (shifted) cohomology
$$HH^{i+1}_R(A):=H^iC^\bullet _R(A)$$
is the {\it Hochschild cohomology} of $A$.

Note that quasi-isomorphic flat minimal $A(\infty)$ algebras have
isomorphic bar constructions (Corollary 2.6), hence isomorphic
Hochschild complexes and Hochschild cohomology.

The Hochschild cohomology $HH_R^\bullet(A)$ is a functor of $R$
which is hard to control because of the presence of infinite
products in the Hochschild complex $C^\bullet_R(A)$. It turns out
that under certain finiteness assumptions on $A$ there is a natural
subcomplex $C^\bullet_{R,c}(A)\subset C^\bullet_R(A)$ whose
cohomology behaves better.

\begin{defi} An $A(\infty)$ algebra $A=(A,(m_1,m_2,...))$ is called
finitely defined if $m_n=0$ for $n>>0$.
\end{defi}

Although the above definition can be made for all $A(\infty)$
algebras (in particular any DG algebra would be a finitely defined
$A(\infty)$ algebra) we think it only makes sense for minimal ones.

For the rest of this section we assume that all $A(\infty)$ algebras
are finitely defined.

\subsection{Definition of $HH_{R,c}^\bullet(A)$}
Recall that the Hochschild complex $C^\bullet_R(A)$ of an
$A(\infty)$ $R$-algebra consists of $R$-modules
$$C_R^p(A)=\prod_{n\geq 1}\Hom ^p_R(A[1]^{\otimes
n},A[1]).$$ Consider the $R$-submodule
$$C_{R,c}^p(A)=\sum_{n\geq 1}\Hom ^p_R(A[1]^{\otimes
n},A[1]).$$ Notice that $C^\bullet _{R,c}(A)$ is actually a
subcomplex of $C_R^\bullet(A)$ since $A$ is finitely defined.

\begin{defi} We call the elements of $C_{R,c}^\bullet(A)$ the
Hochschild cochains with compact supports. The corresponding
cohomology $R$-modules
$$HH_{R,c}^n(A):=H^n(C^\bullet_{R,c}(A))$$
are called the Hochschild cohomology of $A$ with compact supports.
\end{defi}

\subsection{Properties of $HH_{R,c}^\bullet(A)$}

By definition we have the canonical map
$$\iota :HH_{R,c}^\bullet(A)\to HH_R^\bullet(A).$$

\begin{lemma} Assume that $m_n=0$ for $n\neq 2$, i.e. $A$ is just a
graded associative $R$-algebra. Then the map $\iota$ is injective.
\end{lemma}

\begin{proof} Suppose that $d=(d_1,d_2,...)\in C^\bullet_R(A)$ is
a coderivation such that $[m_A,d]=e=(e_1,...,e_n,0,0,...)\in
C^\bullet_{R,c}(A)$. Consider the coderivation $d_{\leq
n-1}:=(d_1,...,d_{n-1},0,0,...)\in C^\bullet_{R,c}(A)$. Then
$[m_A,d_{\leq n-1}]=e$ (because $m_n=0$ for $n\neq 2$), i.e. $e$ is
also a coboundary in the complex $C^\bullet_{R,c}(A)$.
\end{proof}

\begin{prop} Assume  that $A$ is a finite $R$-module. Let
$R\to Q$ be a homomorphism of commutative rings and put
$A_Q=A\otimes _RQ$. Then

a) $C^\bullet_{Q,c}(A_Q)=C^\bullet_{R,c}(A)\otimes _RQ$;

b) If $Q$ is a flat $R$-module, then
$HH_{Q,c}^\bullet(A_Q)=HH_{R,c}^\bullet (A)\otimes _RQ.$
\end{prop}

\begin{proof} Clearly a)$\Rightarrow$b). To prove a) notice the
isomorphism of $Q$-modules
$$\Hom _Q(A_Q^{\otimes _Qn},A_Q)=\Hom _R(A^{\otimes _Rn},A_Q)=\Hom
_R(A^{\otimes _Rn},A)\otimes _RQ$$ (since $A^{\otimes _Rn}$ is a
finite projective $R$-module).
\end{proof}

\begin{remark} In particular, if $A$ is a finite $R$-module then for each $n$
we obtain a quasi-coherent sheaf $\cH \cH _c^n(A)$ on $SpecR$ which
is the localization of the $R$-module $HH_{R,c}^n(A)$.
\end{remark}

\begin{prop} Assume that the ring $R$ is noetherian, $A$ is a finite
$R$-module, and $m_n=0$ for $n\neq 2$ (i.e. $A$ is just a graded
associative $R$-algebra). Also assume that each $R$-module
$HH_{R,c}^n(A)$ is projective. Let $R\to Q$ be a homomoprhism of
commutative rings and put $A_Q=A\otimes _RQ$. Then
$$HH_{Q,c}^n(A_Q)=HH_{R,c}^n(A)\otimes _RQ.$$
\end{prop}

\begin{proof} Since $A$ is just a graded associative algebra, the complex
$C_{R,c}^\bullet(A)$ is a direct sum of complexes
$$C_{R,c}^\bullet(A)=\bigoplus _{i\in \bbZ}C^\bullet _i(A),$$
where $C_i^j(A)=\Hom _R^{i+j}(A^{\otimes j},A).$ Similarly
$$C_{Q,c}^\bullet(A_Q)=\bigoplus _{i\in \bbZ}C^\bullet _i(A_Q).$$
By Proposition 3.4 $C_{Q,c}^\bullet(A_Q)=C_{R,c}^\bullet(A)\otimes
_RQ$ and this isomorphism preserves the decomposition $C^\bullet
=\oplus C^\bullet _i$. So it suffices to prove that for each $i\in
\bbZ$ the complex of $R$-modules $C^\bullet _i(A)$ is homotopy
equivalent to its cohomology $\oplus _nH^n(C^\bullet _i(A))[-n]$
(with the trivial differential). We need a lemma.

\begin{lemma} Let $R$ be a commutative noetherian ring and let
$$K^\bullet
:=...\stackrel{d^{n-1}}{\to}K^n\stackrel{d^n}{\to}K^{n+1}...$$ be a
bounded below complex of finite projective $R$-modules such that
each $R$-module $H^n(K^\bullet)$ is also projective. Then for each
$n$ the $R$-module $\im d^n$ is projective.
\end{lemma}

\begin{proof} Being a projective module is a local property, so we
may and will assume that $R$ is a local noetherian ring. We also may
assume that $K^n=0$ for $n<0$.

Recall the Auslander-Buchsbaum formula: if $M$ is a finite
$R$-module of finite projective dimension $\pd M$ then
$$\pd M+\depth M=\depth R.$$
In particular $\pd M\leq \depth R$.

First we claim that $\pd \im d^n<\infty$ for any $n$. Indeed,
consider the complex
$$0\to K^0\stackrel{d^0}{\to}K^1\stackrel{d^1}{\to}
...\stackrel{d^{n-1}}{\to}K^n\to \im d^n\to 0.$$ This may not be a
projective resolution of $\im d^n$ (since the complex $K^\bullet$
may not be exact), but we can easily make it into one:
$$0\to H^0(K^\bullet)\to K^0\oplus H^1(K^\bullet)\to K^1\oplus
H^2(K^\bullet)\to ...\to K^{n-1}\oplus H^n(K^\bullet)\to K^n\to \im
d^n\to 0$$ where the differential $H^i(K^\bullet)\to K^i$ is any
splitting of the projection $\Ker d^i\to H^i (K^\bullet)$. Thus we
have $\pd \im d^n\leq n$ hence in particular $\pd \im d^n \leq
\depth R$.

But we claim that in fact $\pd \im d^n=0$. The proof is similar.
Indeed, put $\delta =\depth R$ and consider the complex
$$0\to \im d^n\hookrightarrow K^{n+1}\stackrel{d^{n+1}}{\to}
...\to K^{n+\delta}\stackrel{d^{n+\delta}}{\to}\im d^{n+\delta}\to
0.$$ Again we can turn it into an exact complex
$$0\to \im d^n\oplus H^{n+1}(K^\bullet)\to K^{n+1}\oplus
H^{n+2}(K^\bullet)\to ...\to K^{n+\delta}\to \im d^{n+\delta}\to 0$$
which shows that $\pd (\im d^n\oplus H^{n+1}(K^\bullet))=\pd \im
d^n=0$ (since $\pd \im d^{n+\delta}\leq \delta$). This proves the
lemma.
\end{proof}

The lemma implies that for each $n$ we have
$$K^n\simeq \im d^{n-1}\oplus H^n(K^\bullet)\oplus \im d^n.$$
It follows easily that $K^\bullet$ is homotopy equivalent to its
cohomology $\oplus _nH^n(K^\bullet)[-n]$. Now apply this to
$K^\bullet =C^\bullet _i(A)$.
\end{proof}

\begin{remark} We do not know if Proposition 3.6 remains true
without the assumption that $m_n=0$ for $n\neq 2$.
\end{remark}

The following seemingly trivial example is actually an important
one.

\begin{example} Let $k$ be a field and $R$ be a $k$-algebra. Let $B$
be a finitely defined $A(\infty)$ $k$-algebra such that
$\dim_kB<\infty$. Put $A=B\otimes _kR$. Then for each $n$ we have
$$HH_{R,c}^\bullet(A)=HH_{k,c}^\bullet(B)\otimes _kR$$
and hence in particular the corresponding quasi-coherent $\cO
_{SpecR}$-module $\cH \cH _c^n(A)$ is free. Moreover for any
homomorphism of commutative $k$-algebras $R\to Q$ we have
$$HH_{Q,c}^\bullet(A\otimes _RQ)=HH_{k,c}^\bullet(B)\otimes _kQ= HH_{R,c}^\bullet(A)\otimes
_RQ.$$
In particular, if $x\in SpecR$ is a $k$-point, then
$$HH_{k,c}^\bullet(A_x)=HH_{k,c}^\bullet(B).$$
\end{example}

\subsection{Invariance of $HH_{R,c}(A)$}
Let $A$ and $B$ be two flat minimal $A(\infty)$ $R$-algebras which
are finitely defined. Suppose that $A$ and $B$ are quasi-isomorphic.
It is natural to ask whether $HH^\bullet_{R,c}(A)\simeq
HH^\bullet_{R,c}(B)$? This is so at least when there exist mutually
inverse isomorphisms of the bar constructions $f:\cB A\to \cB B$,
$g:\cB B\to \cB A$, such that $f_n=g_n=0$ for $n>>0$. In particular
this is true if $A$ and $B$ are usual associative graded
$R$-algebras (which are isomorphic).

\section{Kaledin's cohomology class}

We thank the referee for suggesting that the material of this
section be presented in a general context of DG Lie algebras. We do
this in Section 7. (The connection being that the Hochshild complex
of an $A(\infty)$-algebra is naturally a DG Lie algebra.) However,
since we are interested in $A(\infty)$-algebras, we decided to also
present this special case explicitly.

\subsection{} Let $k$ be a field of characteristic zero and $R$ be a
commutative $k$-algebra. For an $R$ module $M$ we denote by $M[[h]]$
the $R[[h]]$-module
$$M[[h]]=\lim _\leftarrow M[h]/h^n=\lim _\leftarrow (M\otimes
_RR[h]/h^n)$$ We call an $R[[h]]$-module $P$ {\it $h$-free complete}
if it is isomorphic to $\bar{P}[[h]],$ where $\bar{P}$ is the
$R$-module $P/h.$

Notice that  $M[[h]]$ is canonically identified with the set of
power series $\Sigma _{i=0}^{\infty}m_ih^i,$ $m_i\in M.$ To get the
analogous identification for an arbitrary $h$-free complete
$R[[h]]$-module one needs to choose a splitting $\bar{P}\to P$ (a
map of $R$-modules).

There is a canonical isomorphism of $R[[h]]$-modules
$$\Hom _{R[[h]]}(M[[h]]\otimes _{R[[h]]}...\otimes
_{R[[h]]}M[[h]],M[[h]])=\{ \Sigma _{i=0}^\infty f_ih^i\vert f_i\in \Hom
_R(M\otimes _R...\otimes _RM,M)\}$$

\subsection{} Let $B$ be an $h$-free complete $R[[h]]$-module which
has a structure of a minimal $A(\infty)$ $R[[h]]$-algebra $(B,m).$
Assume that the minimal $A(\infty)$ $R$-algebra
$(\bar{B},m^{(0)})=B/h$ is flat. Choose a splitting $\bar{B}\to B$
of $R$-modules. Then we can write
$$m=m^{(0)}+m^{(1)}h+m^{(2)}h^2+...$$
for some coderivations $m^{(i)}\in C^1_R(\bar{B})$. Notice that the
Hochshild complex $C^\bullet _{R[[h]]}(B)$ is isomorphic to the
inverse limit of the sequence $\{C^\bullet _{R[h]/h^n}(B/h^n)\}$
where all maps are surjective. In particular
$$HH^\bullet _{R[[h]]}(B)=\lim _{\leftarrow} HH^\bullet
_{R[h]/h^n}(B/h^n).$$

Consider the coderivation
$$\partial _hm=m^{(1)}+2m^{(2)}h+3m^{(3)}h^2+...\in C^1_{R[[h]]}(B).$$
Then
$$[m,\partial _hm]=m\cdot \partial _hm+\partial _hm\cdot
m=\partial _h(m\cdot m)=0,$$ i.e. $\partial _hm$ is a cocycle and
hence defines a cohomology class $[\partial _hm]\in
HH_{R[[h]]}^2(B)$.

\begin{lemma} Let $f:\overline{T}B[1]\to
\overline{T}B[1]$ be a coalgebra automorphism which is the identity
modulo $h$. Put $f(c):=f\cdot c\cdot f^{-1}$ for $c\in C^\bullet
_{R[[h]]}(B).$ Then the cocycles $\partial _h(f(m))$ and $f(\partial
_hm )$ are cohomologous (with respect to the differential
$[f(m),-]$).
\end{lemma}

\begin{proof}  It suffices to show this modulo $h^n$ for
all $n$.

Notice that $f$ has the following canonical decomposition
$$f=...\cdot \exp (g^{(2)}h^2)\cdot \exp (g^{(1)}h)$$
for some coderivations $g^{(1)}, g^{(2)}, ...\in C^0_R(\bar{B})$.
Namely, let $f\equiv \id +f^{(1)}h(\text{mod} h^2)$, where
$f^{(1)}=(f^{(1)}_1,f^{(1)}_2,..)$. Let $g^{(1)}$ be the {\it
coderivation} of degree zero defined by the same sequence, i.e.
$g^{(1)}=(f^{(1)}_1,f^{(1)}_2,...).$ Then the coalgebra
automorphisms $f$ and $\exp (g^{(1)}h)$ are equal modulo $h^2$. Now
replace $f$ by $f\cdot \exp (g^{(1)}h)^{-1}\equiv \id
+f^{(2)}h^2(\text{mod} h^3)$. Let $g^{(2)}$ be the coderivation
$g^{(2)}=(f^{(2)}_1,f^{(2)}_2,...)$, etc.

Fix $n\geq 1$. Then
$$f\equiv \exp (g^{(n-1)}h^{n-1})...\exp (g^{(1)}h)(\text{mod} h^n),$$
and we may and will assume that $f=\exp (gh^i)$ for some
coderivation $g\in C^0_{R}(\bar{B})$. We have
$$\partial _h(f(m))=\partial _hf\cdot m\cdot f^{-1}+f\cdot
\partial _hm\cdot f^{-1}-f\cdot m\cdot f^{-1}\cdot \partial _hf
\cdot f^{-1}.$$ So
$$f\cdot \partial _hm\cdot f^{-1}-\partial _h(f(m))=[f(m),
\partial _hf \cdot f^{-1}].$$ But
$$\partial _hf\cdot f^{-1}=\partial _h(\exp (gh^i))\cdot \exp (-gh^i)
=igh^{i-1},$$ so $\partial _hf \cdot f^{-1}\in C^0_{R[[h]]}(B)$ and
hence $\partial _h(f(m))$ and $f(\partial _hm)$ are cohomologous
modulo $h^n$ with respect to the differential $[f(m),-]$.
\end{proof}

\begin{cor} The class $[\partial _hm]\in HH^2_{R[[h]]}(B)$ is well defined,
i.e. is
independent of the choice of the splitting $\bar{R}\to R$.
\end{cor}

\begin{defi} The class $[\partial _hm]\in HH^2_{R[[h]]}(B)$ is
called the Kaledin class of $B$ and denoted $K_B$.
\end{defi}

\begin{remark} The definition of Kaledin class and the above lemma remain
valid for flat minimal $A(\infty)$ $R[h]/h^{n+1}$-algebras. We
consider the class $K_{B/h^{n+1}}$ of the $A(\infty)$
$R[h]/h^{n+1}$-algebra $B/h^{n+1}$ as an element in
$HH^2_{R[h]/h^n}(B/h^n)$.
\end{remark}

\begin{prop} [Ka] Fix $n\geq 1$. Then the class $K_{B/h^{n+1}}\in
HH^2_{R[h]/h^n}(B/h^n)$ is zero if and only if there exists a
quasi-isomorphism of  $A(\infty)$ $R[h]/h^{n+1}$-algebras
$f:B/h^{n+1}\to \bar{B}[h]/h^{n+1}$ such that $f\equiv (\id
,0,0,...)(\text{mod} h)$.
\end{prop}

\begin{proof} Recall that we identify a quasi-isomorphism of two
minimal flat $A(\infty)$ algebras with an isomorphism of their bar
constructions.

One direction is clear: if $f:B/h^{n+1} \to \bar{B}[h]/h^{n+1}$ is a
quasi-isomorphism which is the identity modulo $h$, then
$K_{B/h^{n+1}}=0$ (since by Lemma 4.1 and Remark 4.4 it corresponds
to $K_{\bar{B}[h]/h^{n+1}}=0$ under $f$).

Suppose $K_{B/h^{n+1}}=0$. By induction on $n$ we know that there
exists a quasi-isomorphism $B/h^n\to \bar{B}[h]/h^n$ which is the
identity modulo $h$. Lift this quasi-isomorphism arbitrarily to an
isomorphism of coalgebras $\overline{T}(B/h^{n+1}[1])\to
\overline{T}(\bar{B}[h]/h^{n+1}[1])$. Then by Lemma 4.1 and Remark
4.4 we may and will assume that
$$m=m_{B/h^{n+1}}= m^{(0)}+m^{(n)}h^n$$ and hence
$K_{B/h^{n+1}}=[nm^{(n)}h^{n-1}]$. Since $K_{B/h^{n+1}}=0$ there
exists a coderivation $g\in C^0_R(\bar{B})$ such that
$$[m,gh^{n-1}]=[m^{(0)},gh^{n-1}]=nm^{(n)}h^{n-1}.$$
Consider the coalgebra automorphism $f=\exp
(n^{-1}gh^n):\overline{T}(\bar{B}[h]/h^{n+1}[1])\to
\overline{T}(\bar{B}[h]/h^{n+1}[1])$. Then $m^{(0)}\cdot f= f\cdot
m$, i.e. $f$ is an isomorphism of the bar constructions $f:\cB
(B/h^{n+1})\to \cB (\bar{B}[h]/h^{n+1})$ and hence is a
quasi-isomorphism from $B/h^{n+1}$ to $\bar{B}[h]/h^{n+1}$ (which is
the identity modulo $h$).
\end{proof}

\begin{cor} In the notation of Proposition 4.5 assume that
$m_{B/h^{n+1}}= m^{(0)}+m^{(n)}h^n$. Then there exists a
quasi-isomorphism of  $A(\infty)$ $R[h]/h^{n+1}$-algebras
$f:B/h^{n+1}\to \bar{B}[h]/h^{n+1}$ such that $f\equiv (\id
,0,0,...)(\text{mod} h)$ if and only if the class $[m^{(n)}]\in
HH^2_R(\bar{B})$ is zero.
\end{cor}

\begin{proof} By Proposition 4.5 there exists such a quasi-isomorphism
$f$ if and only if the class $[nm^{(n)}h^{n-1}]\in
HH^2_{R[h]/h^n}(B/h^n)$ is zero. Clearly, this is equivalent to the
class $[m^{(n)}]\in HH^2_R(\bar{B})$ being zero.
\end{proof}

\section{Deformation to the normal cone}

\subsection{}
Let $k$ be a field of characteristic zero and $R$ be a commutative
$k$-algebra. Let $A=(A,m)$ be a minimal flat $A(\infty)$
$R$-algebra. Consider the $A(\infty)$ $R[h]$-algebra
$\tilde{A}=(A[h],\tilde{m}=(m_2,m_3h,m_4h^2,...))$.

\begin{lemma} The map $\tilde{m}$ indeed defines a structure of an
$A(\infty)$ $R[h]$-algebra on $A[h]$.
\end{lemma}

\begin{proof} The defining equation as in Subsection 2.1 above are homogeneous: after the
substitution of $m_ih^{i-2}$ instead of $m_i$ the equation is
multiplied by $h^{n-3}$.
\end{proof}

Denote by $A(2)$ the $A(\infty)$ $R$-algebra $(A,(m_2,0,0,...))$.

\begin{lemma} We have the following isomoprhisms of $A(\infty)$
$R$-algebras.

a) $\tilde{A}/h\simeq A(2)$,

b) $\tilde{A}/(h-1)\simeq A$.
\end{lemma}

\begin{proof} This is clear.
\end{proof}

\begin{defi} The $A(\infty)$ $R[h]$-algebra $\tilde{A}$ is called
the deformation of $A$ to the normal cone.
\end{defi}

\begin{prop} The $A(\infty)$ $R$-algebras $A$ and $A(2)$ are
quasi-isomorphic if and only if the $A(\infty)$ $R[h]$-algebras
$\tilde{A}$ and $A(2)[h]$ are quasi-isomorphic. That is $A$ is
formal if and only if $\tilde{A}$ is such.
\end{prop}

\begin{proof} Given a quasi-isomorphism $\tilde{f}:\tilde{A}\to
A(2)[h]$ we may reduce it modulo $(h-1)$ to get a quasi-isomorphism
between $A$ and $A(2)$. Vice versa, let $f=(f_1,f_2,...):A\to A(2)$
be a quasi-isomorphism of $A(\infty)$ $R$-algebras. Then
$\tilde{f}=(f_1,f_2h,f_3h^2,...)$ is a quasi-isomorphism between
$\tilde{A}$ and $A(2)[h]$.
\end{proof}

\begin{remark} If $A$ and $A(2)$ are
quasi-isomorphic, then there exists a quasi-isomorphism
$\tilde{f}:\tilde{A}\to A(2)[h]$ which is the identity modulo $h$.
Indeed, the last proof produces an $\tilde{f}$, such that
$\tilde{f}\equiv (f_1,0,0,...)(\text{mod} h)$, where $f_1$ is an
algebra automorphism of $A(2)[h]$. Thus we may take the composition
of $\tilde{f}$ with $(f^{-1}_1,0,0,...)$.
\end{remark}

\begin{defi} The $A(\infty)$ $R$-algebra $A$ is called $n$-formal if
there exists a quasi-isomorphism of $A(\infty)$
$R[h]/h^{n+1}$-algebras $\gamma :\tilde{A}/h^{n+1}\to
A(2)[h]/h^{n+1}$, such that $\gamma \equiv (\id ,0,0,...)(\text{mod}
h)$.
\end{defi}

 Notice that Proposition 4.5 above provides a cohomological criterion for
$n$-formality of $A$:

\begin{cor} a) The $A(\infty)$ $R$-algebra $A$ is $n$-formal if and only if
the Kaledin  class $K_{\tilde{A}/h^{n+1}}\in
HH^2_{R[h]/h^n}(\tilde{A}/h^n)$ is zero.

b) Assume that $m_{\tilde{A}/h^{n+1}}=m_2+m_{n+2}h^n$. Then $A$ is
$n$-formal if and only if $[m_{n+2}]\in HH^2_R(A(2))$ is zero (see
Corollary 4.6).
\end{cor}

The next proposition relates $n$-formality to formality.

\begin{prop} The $A(\infty)$ $R$-algebra $A$ is formal if and only if it
is $n$-formal for all $n\geq 1$.
\end{prop}

\begin{proof} One direction is clear: If $A$ and $A(2)$ are
quasi-isomorphic, then by Proposition 5.4 and Remark 5.5 there
exists a quasi-isomorphism of $A(\infty)$ $R[h]$-algebras
 $\tilde{A}\to A(2)[h]$ which is the identity modulo $h$. It
 remains to reduce this quasi-isomorphism modulo $h^{n+1}$.

 Assume that $A$ is $n$-formal for all $n\geq 1$.
By Proposition 5.4 above it suffices to prove that the $A(\infty)$
$R[h]$-algebras $\tilde{A}$ and $A(2)[h]$ are quasi-isomorphic.

 We will prove by
 induction on $n$ that there exists a sequence of maps $g_2$, $g_3$, ..., where
 $g_i\in \Hom _R^0(A[1]^{\otimes i},A[1])$ so that for each
 $n\geq 2$ the following assertion is true:

\medskip

 \noindent{\bf E(n):} Consider maps $g_i$ as coderivations $g_i=(0,...,0,g_i,0,...)$ of degree zero of the
 coalgebra $\overline{T}\tilde{A}[1]$. Then the coalgebra automorphism
 $$\gamma _n:=\exp (g_nh^{n-1})\cdot ...\cdot \exp
 (g_2h):\overline{T}\tilde{A}[1]\to \overline{T}\tilde{A}[1]$$
 when reduced modulo $h^{n}$ becomes a quasi-isomorphism between
 $\tilde{A}/h^{n}$ and $A(2)[h]/h^{n}$.

\medskip

Then the infinite composition $\tilde{f}:=...\exp (g_3h^{2})\exp
 (g_2h)$ is the required quasi-isomorphism between $\tilde{A}$ and
 $A(2)[h]$.

In order to prove the existence of the $g_i$'s it is convenient to
introduce the following $k^*$-action on the $R$-module
$\overline{T}\tilde{A}[1]$. For $\lambda \in k^*$ put
$$\lambda \star x:=\lambda ^i x ,\quad \text{if $x\in
(A[1])^{\otimes i}$},\quad \text{and $\lambda \star h=\lambda h$}.$$
Notice that both $m_2$ and $\tilde{m}$ are maps of degree $-1$ with
respect to this action.

Now assume that we found $g_2$, ...,$g_n$ so that {\bf E(n)} holds.
Then
$$\gamma _n \cdot \tilde{m} \cdot \gamma _n^{-1}\equiv
m_2+m_{n}^\prime h^{n}(\text{mod} h^{n+1})$$ for some coderivation
$m_{n}^\prime \in C^1_R(A)$.
 Notice that the map $\gamma _n$ is of degree zero with respect to
the $k^*$-action. Hence the coderivation $\gamma _n \cdot \tilde{m}
\cdot \gamma _n^{-1}$ is again of  degree $-1$. This forces the
coderivation $m_{n}^\prime $ to be defined by a single map in $\Hom
^1_R(A[1]^{\otimes n+2},A[1])$. Since $A$ is $n$-formal, by
Corollary 4.6 the class  $[m_n^\prime]$ is zero in $HH^2_R(A(2))$.
So there exists a coderivation $g_{n+1}\in C^0_R(A)$ such that
$[m_2,g_{n+1}]=m_{n}^\prime$. It is clear that we can choose
$g_{n+1}$ to be defined by a single map $g_{n+1}\in \Hom
^0_R(A[1]^{\otimes n+1},A[1])$. Then the coalgebra isomorphism
$$\gamma _{n+1}:=\exp (g_{n+1}h^n)\cdot \gamma _n:\overline{T}\tilde{A}[1]\to \overline{T}\tilde{A}[1]$$
induces a quasi-isomorphism between $\tilde {A}/h^{n+1}$ and
$A(2)[h]/h^{n+1}$. This completes the induction step and proves the
proposition.
\end{proof}

\subsection{}
Notice that for each $n\geq 1$ the $A(\infty)$ algebra
$\tilde{A}/h^{n}$ is finitely defined. Thus the Hochshild cohomology
with compact supports $HH_{R[h]/h^n,c}^\bullet (\tilde{A}/h^n)$ is
defined. Moreover the Kaledin class $K_{\tilde{A}/h^{n+1}}$
obviously belongs to the image of $HH_{R[h]/h^n,c}^\bullet
(\tilde{A}/h^n)$ in $HH^2_{R[h]/h^n}(\tilde{A}/h^n)$. Therefore it
is useful to notice the following fact.

\begin{lemma} For any $n\geq 1$ the canonical map
$$HH_{R[h]/h^n,c}^\bullet (\tilde{A}/h^n)\to HH_{R[h]/h^n}^\bullet (\tilde{A}/h^n)$$ is injective.
\end{lemma}

\begin{proof} This is easy to see by considering the weights of the
$k^*$-action as in the proof of Proposition 5.8.
\end{proof}

\begin{remark} Thus we may and will consider the obstruction to $n$-formality
of $A$ (i.e. the Kaledin class $K_{\tilde{A}/h^{n+1}}$) as an
element of $HH^2_{R[h]/h^n,c}(\tilde{A}/h^n)$. In particular in
Corollaries 4.6 and 5.7 we can use the Hochschild cohomology with
compact supports.
\end{remark}

\section{Applications}

\subsection{Formality of $A(\infty)$ algebras} Let $k$ be a field of characteristic zero and $R$ be a
commutative $k$-algebra. Let $A=(A,m)$ be a minimal flat $A(\infty)
$ $R$-algebra and $\tilde{A}$ be its deformation to the normal cone.
If $m=(m_2,m_3,...)$ denote as before $A(2):=(A,(m_2,0,0,...))$,
i.e. $A(2)$ is the underlying associative algebra of $A$. We have
$A(2)=\tilde{A}/h$. By definition $A$ is formal if it is
quasi-isomorphic to $A(2)$.

\begin{remark} Let $R\to Q$ be a homomorphism of commutative
$k$-algebras. If $A$ is formal then clearly the $A(\infty)$
$Q$-algebra $A_Q=A\otimes _RQ$ is also formal.
\end{remark}

\begin{prop} Assume that $A$ is a finite $R$-module. Let $R\to Q$ be
a  homomorphism of commutative rings. Put $A_Q=A\otimes _RQ$. Assume
that $Q$ is a faithfully flat $R$-module. Then $A$ is formal if and
only if the $A(\infty)$ $Q$-algebra $A_Q$ is formal.
\end{prop}

\begin{proof} By Proposition 5.8 $A$ (resp. $A_Q$) is
formal if and only if it is $n$-formal for all $n\geq 1$.

Fix $n\geq 1$. Notice that $Q[h]/h^n$ is faithfully flat over
$R[h]/h^n$. By Proposition 3.4 we have
$HH_{Q[h]/h^n,c}^2(\tilde{A_Q}/h^n)=HH_{R[h]/h^n,c}^2(\tilde{A}/h^n)\otimes
_{R[h]/h^n}Q[h]/h^n$. And by faithful flatness the class
$K_{\tilde{A}/h^{n+1}}\in HH^2_{R[h]/h^n,c}(\tilde{A}/h^n)$ is zero
if and only if the class
$K_{\tilde{A_Q}/h^{n+1}}=K_{\tilde{A}/h^{n+1}}\otimes 1\in
HH^2_{Q[h]/h^n,c}(\tilde{A_Q}/h^n)$ is zero . Hence the proposition
follows from Corollary 5.7 a) and Remark 5.10.
\end{proof}

\begin{prop} Assume that $R$ is an integral domain with the generic
point $\eta \in SpecR$. Assume that $A$ is a finite $R$-module and
that the $R$-module $HH_{R,c}^2(A(2))$ is torsion free. If the
$A(\infty)$ $k(\eta)$-algebra $A_\eta$ is formal then $A$ is also
formal. In particular the $A(\infty)$ $k(x)$-algebra $A_x$ is formal
for all points $x\in SpecR$.
\end{prop}

\begin{proof} By Proposition 5.8 it suffices to prove
that $A$ is $n$-formal for all $n\geq 1$. We do it by induction on
$n.$ Fix $n\geq 1$ and assume that $A$ is $(n-1)$-formal. Then we
may and will assume that $m_{\tilde{A}/h^{n+1}}= m_2+m_{n+2}h^n$. By
Corollary 5.7 b) and Remark 5.10 $A$ is $n$-formal if and only if
the class $[m_{n+2}]\in HH^2_{R,c}(A(2))$ is zero. This class
vanishes at the generic point $\eta$ (since $HH_{R,c}^2(A(2))\otimes
_Rk(\eta)=HH^2_{k(\eta),c}(A_{\eta}(2))$ by Proposition 3.4) and
hence vanishes identically, since the $R$-module $HH_{R,c}^2(A(2))$
is torsion free. This completes the induction step and proves the
proposition.
\end{proof}

\begin{prop} Let $R$ be noetherian. Assume that $A$ is a finite $R$-module
and that for each $n$ the $R$-module $HH_{R,c}^n(A(2))$ is
projective. Then the subset
$$F(A):=\{ x\in SpecR\ \vert \ \text{the $A(\infty)$ $k(x)$-algebra $A_x$
is formal}\}$$ is closed under specialization.
\end{prop}

\begin{proof} We may assume that $F(A)$ is not empty.
Choose $\eta \in F(A)$ and consider its closure $\overline
{\eta}=:Spec\bar{R}\subset SpecR$. Then $\bar{R}$ is an integral
domain and $A_{\bar{R}}=A\otimes _R\bar{R}$ is an (flat minimal)
$A(\infty)$ $\bar{R}$-algebra which is a finite $\bar{R}$-module. By
Proposition 3.6 above
$HH_{\bar{R},c}^2(A(2)_{\bar{R}})=HH^2_{R,c}(A(2))\otimes
_R\bar{R}$. This is a projective $\bar{R}$-module, in particular,
torsion free. Hence the assumptions of the previous proposition hold
for $A_{\bar{R}}$ and thus $A_{\bar{R}}$ is formal. So $A_x$ is
formal for all $x\in Spec\bar{R}$.
\end{proof}

\begin{prop} Let $R$ be noetherian and $I\subset R$ be an ideal such that $\cap _nI^n=0$.
Assume that $A$ is a finite $R$-module and for each $n$ the
$R$-module $HH_{R,c}^n(A(2))$ is projective. Assume that the
$A(\infty)$ $R/I^n$-algebra $A_n:=A/(I)^n$ is formal for all $n\geq
1$. Then $A$ is formal.
\end{prop}

\begin{proof} The proof is similar to the proof of Proposition 6.3
Namely we prove by induction on $n$ that $A$ is $n$-formal. Fix
$n\geq 1$ and assume that $A$ is $n-1$-formal. Then we may assume
that $m_{\tilde{A}/h^{n+1}}= m_2+m_{n+2}h^n$. By Corollary 5.7 b)
and Remark 5.10 $A$ is $n$-formal if and only if the class
$[m_{n+2}]\in HH^2_{R,c}(A(2))$ is zero. By Proposition 3.6 we have
$$HH_{R,c}^2(A(2))\otimes _RR/I^l=HH^2_{R/I^l,c}(A(2)/I^l)$$
and by our assumption the class $[m_{n+2}]\otimes 1\in
HH^2_{R/I^l,c}(A(2)/I^l)$ is zero.
 Therefore the class $[m_{n+2}]=0$, because $\cap _lI^l=0$ and the $R$-module $HH^2_{R,c}(A(2))$ is projective.
 This completes the induction step and
proves the proposition.
\end{proof}

\begin{prop} Assume that $R$ is noetherian and has the trivial radical
(i.e. the intersection of maximal ideals of $R$ is zero). Assume
that $A$ is a finite $R$-module. Assume that for each $n$ the
$R$-module $HH_{R,c}^n(A(2))$ is projective. If $A_x$ is formal for
all closed points $x\in SpecR$ then $A$ is formal (and hence $A_y$
is formal for all points $y\in SpecR$).
\end{prop}

\begin{proof}
Again we use Proposition 5.8: it suffices to prove that $A$ is
$n$-formal for all $n\geq 1$. Fix $n\geq 1$ and assume that $A$ is
$n-1$-formal. Then we may assume that $m_{\tilde{A}/h^{n+1}}=
m_2+m_{n+2}h^n$. By Corollary 5.7 b) and Remark 5.10 $A$ is
$n$-formal if and only if the class $[m_{n+2}]\in HH^2_{R,c}(A(2))$
is zero. Let $J\subset R$ be a maximal ideal. By Proposition 3.6 we
have
$$HH_{R,c}^2(A(2))\otimes _RR/J=HH^2_{R/J,c}(A(2)/J)$$
and by our assumption the class $[m_{n+2}]\otimes 1\in
HH^2_{R/J,c}(A(2)/J)$ is zero. Therefore the class $[m_{n+2}]=0$,
because the radical of $R$ is trivial and $HH^2_{R,c}(A(2))$ is a
projective $R$-module. This completes the induction step and proves
the proposition.
\end{proof}

\begin{remark} Assume that there exists an associative graded
$k$-algebra $B$ such that the $A(2)=B\otimes _kR$ and $\dim _kB
<\infty.$ Then we may consider $A$ as an $R$-family of
$A(\infty)$-structures which extend the same associative algebra
structure on $B$. In this case for each $n$ the $R$-module
$HH_{R,c}^n(A(2))$ is free and the conclusions of Proposition 6.4,
6.5, 6.6 hold without the assumption of $R$ being noetherian
(Example 3.9).
\end{remark}

\subsection{Formality of DG algebras}
All the results of this section can be formulated in the language of
DG algebras rather than $A(\infty)$ algebras. Namely, assume again
that $k$ is a field of characteristic zero and $R$ be a commutative
$k$-algebra. Let $\cA$ be flat DG $R$-algebra, i.e. each cohomology
$R$-module $H^i(\cA)$ is {\it projective.} Then by Theorem 2.2 it
has a minimal $A(\infty)$ model $A,$ which is unique up to a
quasi-isomorphism (Corollary 2.6). It comes with an $A(\infty)$
quasi-isomorphism $A\to \cA.$  By Corollary 2.9 $\cA$ is formal (as
a DG algebra) if and only if $A$ is formal (as an $A(\infty)$
algebra).

We would like to study extended DG algebras $\cA\otimes _RQ,$ for
various (commutative) algebra homomorphisms $R\to Q.$ In particular
we would like to study the fibers $\cA _x$ of $\cA$ at various
points of $x\in Spec R.$ To do that we should first replace the DG
algebra $\cA$ by a quasi-isomorphic one which is cofibrant.

\begin{lemma} Let $\cC$ be a cofibrant DG $R$-algebra. Then $\cC$ is
cofibrant as a complex of $R$-modules.
\end{lemma}

\begin{proof} This follows from [Sch-Sh], Theorem 4.1(3).
Alternatively, it is easy to see directly if $\cC$ is semi-free
([Dr]).
\end{proof}

So from now on we assume that the flat DG algebra $\cA$ is
cofibrant. The the $A(\infty)$ quasi-isomorphism $A\to \cA$ remains
a quasi-isomorphism after any extension of scalars.

\begin{cor} Let $\cA$
be DG $R$-algebra such that the total cohomology $R$-module
$H^\bullet(\cA)$ is  projective of finite rank. Let $R\to Q$ be a
homomorphism of commutative rings. Assume that $Q$ is a faithfully
flat $R$-module. Then $\cA$ is formal if and only if the DG
$Q$-algebra $\cA \otimes _RQ$ is formal.
\end{cor}

\begin{proof} Let $A$ be a minimal $A(\infty)$ $R$-algebra with a
quasi-isomorphism of $A(\infty)$ $R$-algebras $f:A \to \cA$. Then
$f\otimes \id :A\otimes _RQ\to \cA \otimes _RQ$ is also a
quasi-isomorphism. So the corollary follows from Proposition 6.2
\end{proof}

\begin{cor} Let $\cA$
be DG $R$-algebra such that total cohomology $R$-module
$H^\bullet(\cA)$ is  projective of finite rank and $\cA$ is
cofibrant as a complex of $R$-modules. We consider the cohomology
$H^\bullet (\cA)$ as an $A(\infty)$ algebra with $m_i=0$ for $i\neq
2.$

a)  Assume that $R$ is an integral domain with the generic point
$\eta \in SpecR$. Assume that the $R$-module $HH_{R,c}^2(H^\bullet
(\cA))$ is torsion free. If the DG $k(\eta)$-algebra $\cA_\eta$ is
formal then the DG $R$-algebra $\cA$ is also formal. In particular,
$\cA_x$ is formal for all points $x\in SpecR$.

b) Let $R$ be noetherian. Assume that for each $n$ the $R$-module
$HH_{R,c}^n(H^\bullet (\cA))$ is projective. Then the subset
$$F(\cA):=\{ x\in SpecR\ \vert \ \text{the DG $k(x)$-algebra $\cA_x$
is formal}\}$$ is closed under specialization.

c) Let $R$ be noetherian and $I\subset R$ be an ideal such that
$\cap _nI^n=0$. Assume that for each $n$ the $R$-module
$HH_{R,c}^n(H^\bullet (\cA))$ is projective. Assume that the DG
$R/I^n$-algebra $\cA \otimes _RR/I^n=\cA/(I)^n$ is formal for all
$n\geq 1$. Then $\cA$ is formal.

d) Assume that $R$ is noetherian and has the trivial radical (i.e.
the intersection of maximal ideals of $R$ is zero).  Assume that for
each $n$ the $R$-module $HH_{R,c}^n(H^\bullet (\cA))$ is projective.
If $\cA_x$ is formal for all closed points $x\in SpecR$ then $\cA$
is formal (and hence $\cA_y$ is formal for all points $y\in SpecR$).
\end{cor}

\begin{proof} This follows from
Propositions 6.3, 6.4, 6.5, 6.6 above. Indeed, if $A\to \cA$ is a
minimal flat $A(\infty)$ model for $\cA,$ then $H^\bullet
(\cA)=A(2)$ and for any  homomorphism $R\to Q$ of commutative
algebras the DG $Q$-algebra $\cA \otimes _RQ$ is DG formal if and
and only if the $A(\infty)$ $Q$-algebra $A\otimes _RQ$ is formal.
\end{proof}

\begin{remark} Let $\cA$ be as in the last corollary. Assume that
there exists an associative $k$-algebra $B$ such that $H^\bullet
(\cA)=B\otimes _kR$. Then we may consider $\cA$ as an $R$-family of
DG algebras with the "same" cohomology algebra. In this case for
each $n$ the $R$-module $HH^n_{R,c}(H^\bullet (\cA))$ is free and
the conclusions in parts b),c),d) of the corollary hold without the
assumption of $R$ being noetherian (Remark 6.7).
\end{remark}

\section{Kaledin cohomology class for DG algebras}

\subsection{DG Lie algebras}

Let $k$ be a field of characteristic zero, $R$ be a commutative
$k$-algebra and $L=\oplus L^i$ be a graded $R$-module. Assume that
there is given an $R$-linear map $[\ , \ ]:L\otimes _RL\to L$ which
is homogeneous of degree zero and satisfies the following relations
$$[\alpha ,\beta ]+(-1)^{\bar{\alpha}\bar{\beta}}[\beta
,\alpha]=0,$$
$$(-1)^{\bar{\gamma}\bar{\alpha}}[\alpha,[\beta ,\gamma]]+
(-1)^{\bar{\alpha}\bar{\beta}}[\beta,[\gamma ,\alpha]]+
(-1)^{\bar{\beta}\bar{\gamma}}[\gamma,[\alpha ,\beta]]=0,$$ where
$\bar{x}$ denotes the degree of a homogeneous element $x\in L$. Then
$L$ is called a graded Lie $R$-algebra.

A homogeneous $R$-linear map $d:L \to L$ of degree $l$ is called a
derivation if
$$d([\beta ,\gamma])=[d\beta ,\gamma ]+(-1)^{\bar{\beta}l}[\beta
,d\gamma].$$ Homogeneous $R$-linear derivations of $L$ form a graded
Lie algebra
$$\Der _R(L)=\Der (L)=\oplus \Der ^i(L).$$

We have a natural homomorphism of graded algebras
$$\ad :L\to \Der (L), \quad \ad _{\alpha }(-):=[\alpha ,-].$$

\begin{defi} A DG Lie algebra is a pair (L,d), where $L$ is a graded
Lie algebra and $d\in \Der ^1(L)$ is such that $d^2=0$.
\end{defi}

Notice that the cohomology of a DG Lie algebra is naturally a graded
Lie algebra.

\subsection{Gauge group}
Let $\frak{g}$ be an graded Lie $R$-algebra. Consider the graded Lie
$R[[h]]$-algebra
$$\frak{g}[[h]]:=\oplus _i\frak{g}^i[[h]],$$
where $\frak{g}^i[[h]]$ consists of power series $\alpha _0+\alpha
_1h+\alpha _2h^2+...,$ $\alpha _n\in \frak{g}^i$ with the bracket
induces by $[\alpha h^n,\beta h^m]=[\alpha ,\beta ]h^{n+m}$.

Clearly $\frak{g}^i[[h]]=\lim_{\leftarrow}\frak{g}^i[[h]]/h^n$ for
each $i$. In particular, the Lie subalgebra $h\frak{g}^0[[h]]\subset
\frak{g}[[h]]$ is the inverse limit of nilpotent Lie algebras
$$\frak{g}^0_n:=h\frak{g}^0[[h]]/h^{n+1}.$$

Let $G_n$ be the group of $R[h]$-linear automorphisms of the graded
Lie algebra $\frak{g}[[h]]/h^{n+1}$ generated by operators $\exp
^{\ad _{\alpha}},$ $\alpha \in h\frak{g}^0[[h]]/h^{n+1}$ which act
by the formula
$$\exp (\ad _{\alpha})(\beta )=\beta +[\alpha
,\beta] +\frac{1}{2!}[\alpha,[\alpha,\beta]]+...$$

Notice that by the Campbell-Hausdorff formula every element of $G_n$
is equal to $\exp ^{\ad _{\alpha}},$ for some $\alpha \in
h\frak{g}^0[[h]]/h^{n+1}$.

There are natural surjective group homomorphisms $G_{n+1}\to G_n$
and we denote
$$G=G(\frak{g}):=\lim _{\leftarrow}G_n.$$
The group $G$ is called the gauge group of $\frak{g}$. It acts
naturally by $R[[h]]$-linear automorphisms of the graded Lie algebra
$\frak{g}[[h]]$ by the adjoint action. This action is by definition
faithful. This induces the action of $G$ on the graded Lie algebra
$\Der (\frak{g}[[h]])$. In particular, if $(\frak{g}[[h]],d)$ is a
DG Lie algebra and $g\in G,$ then $(\frak{g}[[h]],g(d))$ is also
such.

\subsection{Kaledin class} Let $(\frak{g},d)$ be a DG Lie
$R$-algebra. Consider the DG Lie $R[[h]]$-algebra
$(\frak{g}[[h]],d)$. Let $\pi=\pi _1h+\pi _2h^2+... \in
h\frak{g}^1[[h]]$ be a solution of the Maurer-Cartan equation
$$d\pi +\frac{1}{2}[\pi,\pi]=0.$$
In other words the derivation $d_{\pi}:=d+[\pi ,-]$ satisfies $d_\pi
^2=0$. Consider the element
$$\partial _h(d_{\pi})=\partial _h(\pi)=\pi _1+2\pi _2h+3\pi
_3h^2+...\in \frak{g}^1[[h]].$$ We have
$$0=\partial _h(d_{\pi}^2)=\partial _h(d_{\pi})\cdot
d_{\pi}+d_{\pi}\cdot \partial _h(d_\pi)=[d_\pi ,\partial
_h(d_{\pi})].$$ Thus $\partial _h(d_{\pi })$ is a 1-cocycle in the
DG Lie algebra $(\frak{g}[[h]],d_\pi)$.

\begin{defi} We call the corresponding cohomology class $[\partial
_h(d_\pi)]\in H^1(\frak{g}[[h]],d_\pi)$ the Kaledin class (of
$\pi$).
\end{defi}

\begin{prop} a) The Kaledin class $[\partial
_h(d_\pi)]\in H^1(\frak{g}[[h]],d_\pi)$ is gauge invariant. That is
for $g\in G$ the classes $[g(\partial _h(d_{\pi}))],\ [\partial
_h(g(d_\pi))]\in H^1(\frak{g}[[h]],g(d_{\pi}))$ are equal.

b) Moreover, the class $[\partial _h\pi ]=0$ if and only if $\pi $
is gauge equivalent to zero, i.e. there exists $g\in G$ such that
$g(d_{\pi})=d.$
\end{prop}

\begin{proof} a).  Since
$$H^\bullet (\frak{g})=\lim _{\leftarrow}H^\bullet
(\frak{g}[[h]]/h^n)$$ it suffices to prove that the two classes are
congruent modulo $h^{n+1}$ for all $n\geq 0$. So fix $n\geq 0$ and
$g\in G$. Since we work modulo $h^{n+1}$ we may and will assume that
$g\in G_n$.

\begin{lemma} There exist $\xi _1,...,\xi _n\in \frak{g}^0$ such
that
$$g=\exp(\xi _nh^n)\exp(\xi _{n-1}h^{n-1})...\exp(\xi _1h).$$
\end{lemma}

\begin{proof} By induction on $n$ we assume that the statement of
the lemma holds for the image of $g$ in $G_{n-1}$. Thus there exist
$\xi _1, ..., \xi _{n-1} \in \frak{g}^0$ so that
$$\bar{g}:=\exp (-\xi _1h)...\exp (-\xi _{n-1}h^{n-1})g$$
lies in the kernel of the projection $G_n\to G_{n-1}$. Let $\eta =
\eta _1h+...+\eta _nh^n\in h\frak{g}^0[[h]]/h^{n+1}$ be such that
$\bar{g}=\exp(\eta)$. Since the image of $\bar{g}$ under the
projection $G_n\to G_1$ is trivial we conclude that $\eta _1$ is in
the center of the graded Lie algebra $\frak{g}$. Hence we may and
will assume that $\eta _1=0$. Similarly, considering the trivial
image of $\bar{g}$ under the projection $G_n\to G_2$ we may and will
assume that $\eta _2=0$, etc. So $\bar{g}=\exp (\eta _nh^n)$ and we
can take $\xi _n=\eta _n$. This proves the lemma.
\end{proof}

Using the lemma we may and will assume that $g=\exp(\xi h^i)$ for
some $\xi \in \frak{g}^0$, $i>0$. By definition $g(d_{\pi})=g\cdot
d_{\pi}\cdot g^{-1},$ hence
$$\partial _h(g(d_{\pi}))=\partial _hg\cdot d_{\pi}\cdot g^{-1}+
g\cdot \partial _h(d_{\pi})\cdot g^{-1}-g\cdot d_{\pi}\cdot g^{-1}
\cdot\partial _hg\cdot g^{-1}.$$ So
$$g(\partial _h(d_{\pi}))-\partial
_h(g(d_{\pi}))=[g(d_{\pi}),\partial g\cdot g^{-1}].$$ But
$$\partial _hg\cdot g^{-1}=\partial _h(\exp (\xi h^i))\cdot \exp (-\xi h^i)
=i\xi h^{i-1}.$$ This proves a).

b). If $g(d_{\pi})=d$ for some $g\in G$, then $\partial
_h(g(d_{\pi}))=0$ and hence by part a) also $[\partial _h(\pi)]=0$.

Vice versa, suppose that $[\partial _h(\pi)]=0$. Let $\pi =\pi
_1h+\pi _2h^2+...$. Then in particular $0=[\pi _1]\in
H^1(\frak{g},d)$. So there exists $\xi_1\in \frak{g}^0$ such that
$d(\xi _1)=\pi _1$. Put $g_1:=\exp (\xi _1h)\in G.$ Then
$$g_1(d_{\pi})\cong d(\text{mod} h^2).$$
By induction we may assume that we found $\xi _1,...,\xi _{n-1}\in
\frak{g}^0$ so that
$$g_{n-1}...g_1(d_{\pi})\cong d(\text{mod} h^n),$$
where $g_i=\exp (\xi _ih^i).$ Then by part a) we may assume that
$\pi _1=...=\pi _{n-1}=0.$ So by our assumption we have in
particular $0=[n\pi _nh^{n-1}]\in H^1(\frak{g}[[h]]/h^n ,d_{\pi }).$
This is equivalent to saying that $0=[n\pi _n]\in H^1(\frak{g},d).$
Let $\xi _n\in \frak{g}^0$ be such that $d(\xi _n)=[\pi _n]$ (recall
that $\bbQ \subset R$) and put $g_n:=\exp (\xi _n)$. Then
$$g_n(d_{\pi})\cong d(\text{mod} h^{n+1}).$$
This completes our induction step. Put $g:=...g_3g_2g_1\in G$. Then
$$g(d_{\pi})=d.$$
\end{proof}

If we consider the DG Lie $R[[h]]$-algebra $(\frak{g}[[h]],d_{\pi})$
as a deformation of the DG Lie $R$-algebra $(\frak{g}, d)$, then
Proposition 7.3 above asserts that this deformation is trivial if
and only if the Kaledin class $[\partial _h(d_{\pi})]\in
H^1(\frak{g}[[h]],d_{\pi})$ is zero.

All the above can be repeated for DG Lie $R[[h]]/h^n$-algebras
$(\frak{g}[[h]]/h^n,d_{\pi})$. In particular we obtain the following
corollary.

\begin{cor} a) The Kaledin class $[\partial _h(d_{\pi })]\in
H^1(\frak{g}[[h]]/h^{n+1},d_{\pi})$ is gauge invariant, i.e. for
$g\in G_n$ the classes $[g(\partial _h(d_{\pi}))],\ [\partial
_h(g(d_\pi))]\in H^1(\frak{g}[[h]]/h^{n+1},g(d_{\pi}))$ are equal.

b) Moreover, the class $[\partial _h\pi ]=0$ if and only if $\pi $
is gauge equivalent to zero, i.e. there exists $g\in G_n$ such that
$g(d_{\pi})=d.$
\end{cor}

\begin{proof} Same as that of Proposition 7.3.
\end{proof}


\begin{thebibliography}{99}

\bibitem[Dr]{Dr}V.~Drinfeld, DG quotients of DG categories, J. Algebra 272 (2004), no. 2, 643-691.  ArXiv:math/0210114.

\bibitem[Ha-St]{Ha-St}S.~Halperin, J.~Stasheff, Obstructions to
homotopy equivalences, {\it Advances in Math.} No. 32, 233-279
(1979).

\bibitem[Hi]{Hi}V.~Hinich, Tamarkin's proof of Kontsevich formality
theorem, Forum Math. 15 (2003), no. 4, 591-614.  ArXiv:math/0003052.


\bibitem[Ka]{Ka}D.~Kaledin, Some remarks on formality in families,
Mosc. Math. J. 7 (2007), no. 4, 643-652. ArXiv:math/0509699.

\bibitem[Kad1]{Kad}T.~V.~Kadeishvili, Algebraic structure in the
homolo gies of an $A(\infty )$ algebra, (in Russian), {\it Bulletin
of Acad. of Sci. of Georgian SSR}, 108, No. 2, (1982).

\bibitem[Kad2]{Kad}T.~V.~Kadeishvili, The functor D for the categry of $A(\infty)$-algebras   (in Russian), {\it
Soobshch. Acad. Nauk Gruzin. SSR}, 125 (1987), 273-276.

\bibitem[Ke]{Ke}B.~Keller, Introduction to $A$-infinity algebras and
modules, Homology Homotopy Appl. 3 (2001), no. 1, 1-35.
ArXiv:math/9910179.

\bibitem[Le-Ha]{Le-Ha}K.~Lef\`{e}vre-Hasegava, Th\`{e}se de
Doctorat, Universit\'{e} Paris 7 (2003)

\bibitem[Sch-Sh]{Sch-Sh}S.~Schwede, B.~Shipley, Algebras and modules in monoidal model categories,
Proc. London Math. Soc. (3) 80 (2000) 491-511.

\bibitem[Su]{Su}D.~Sullivan, Infinitesimal computations in topology, {\it Publications de IHES}
No. 47.


\end{thebibliography}
\end{document}